\documentclass[3p,times]{elsarticle}
\usepackage{amssymb}
\usepackage{amsxtra}
\usepackage{amsmath}
\usepackage{amsfonts}
\usepackage{CJK}
\usepackage[titletoc]{appendix}
\usepackage{graphicx}

\numberwithin{equation}{section}

\newcommand{\eqa}{\begin{eqnarray}}
\newcommand{\eeqa}{\end{eqnarray}}
\newcommand{\beq}{\begin{equation}}
\newcommand{\eeq}{\end{equation}}
\newcommand{\nn}{\nonumber}

\allowdisplaybreaks

\begin{document}
\begin{frontmatter}
\title{Riemann-Hilbert method and soliton solutions in the system of two-component Hirota equations}
\author[ZL]{Fang Fang\corref{cor1}}
\author[ZL]{Beibei Hu\corref{cor1}}
\author[ZL]{Ling Zhang}
\author[ZN]{Ning Zhang}

\cortext[cor1]{Corresponding authors. E-mail addresses:fangfang7679@163.com (F. Fang), hu\_chzu@shu.edu.cn(B.-B. Hu).}
\address[ZL]{School of Mathematics and Finance, Chuzhou University, Anhui 239000, China}
\address[ZN]{Department of Basical Courses, Shandong University of Science and Technology, Taian 271019, China}

\pagestyle{plain}
\setcounter{page}{1}
\begin{abstract}
In this letter we examine the two-component Hirota (TH) equations which describes the pulse propagation in a coupled fiber with higher-order dispersion and self-steepening. As the TH equations is a complete integrable system, which admits a $3\times 3$ Ablowitz-Kaup-Newell-Segu(AKNS)-type Lax pair, we obtain the general N-soliton solutions of the TH equations via the Riemann-Hilbert(RH) method when the jump matrix of a specific RH problem is a $3\times3$ unit matrix.  As an example, the expression of one- and two-soliton are displayed explicitly.
\end{abstract}
\end{frontmatter}

\section{Introduction}

Soliton theory is a crucial research content in nonlinear science. Due to the significant application of the soliton theory in mathematics and physics, its research has received universal attention by physicists and mathematicians. such as, soliton theory provides a series of methods for solving integrable linear and nonlinear evolution partial differential equations (PDEs in brief) in mathematics and the solitons are often used to describe solitary waves with elastic scattering properties in physics. In many disciplines, there are problems related to soliton theory, which makes it paramount to establish soliton equation model and construct its analytical solution, especially soliton solution. With the development of the soliton theory, increasing methods for solving soliton equations have emerged. such as the inverse scattering transform (IST in brief) \cite{GGMM1967,MJA1981,RB1984,PAC1989}, the Hirota method \cite{HR1971,MLY2016}, the B\"{a}cklund transformation method \cite{PJO1986}, the Darboux transformation (DT in brief) method \cite{VBM1991,LCZ2013,HJS2013,WLH2019} and other methods \cite{LSY1990,LJB2000,FEG2003,CK2018}.
Recently, a new powerful method, the Riemann-Hilbert (RH in brief) method has been developed to the examine of N-soliton solutions \cite{YJK2010}. Through this method, the N-soliton solution for an ocean of integrable systems are obtained \cite{GBL2012,WZ2016,AB2017,ZYS2017,MWX2018,HJ2018,ALP2018}. In particular, the RH method is an effective way to working the initial-boundary value problem of the integrable nonlinear evolution PDEs \cite{Hu1,Hu2,ZQZ2018}.

It is well know that the two-component Hirota (TH in brief) equations can be effect explains pulse propagation in single mode fibers which reads \cite{BSG2001}
\beq\begin{array}{l}
k_1q_{1t}+2A_1k_1q_{1xx}+4k_1^3A_1(|q_1|^2+|q_2|^2)q_1
+i\epsilon[k_1q_{1xxx}+3ik_1^3(|q_1|^2+|q_2|^2)q_{1x}+3ik_1^3q_1(q_1^*q_{1x}+q_2^*q_{2x})]=0,\\
k_1q_{2t}+2A_1k_1q_{2xx}+4k_1^3A_1(|q_1|^2+|q_2|^2)q_2
+i\epsilon[k_1q_{2xxx}+3ik_1^3(|q_1|^2+|q_2|^2)q_{2x}+3ik_1^3q_2(q_1^*q_{1x}+q_2^*q_{2x})]=0,
\end{array}\label{1.1}\eeq
where $q_j(x,t)$ is the complex smooth envelops, and $\epsilon$ represent the strength of high-order effects.
Indeed, when $k_1=1,A_1=-\frac{i}{2}$, that the above system \eqref{1.1} is the bright soliton version of the TH equations, which the Lax pair and the IST method were reported in \cite{TSR1992}, and N-soliton solutions has been discussed via Hirota bilinear form \cite{RR1996}, and rogue wave solutions were obtained by using of DT \cite{PK1997}, and the bright soliton solitons are discussed by RH formulation in \cite{WDS2014}.

On the other hand, when $k_1=i,A_1=\frac{i}{2}$, that the above system \eqref{1.1} is the dark soliton version of the TH equations, which the Painlev\'{e} analysis, the dark soliton solutions and the Lax pair for the N-coupled Hirota equations have been studied \cite{BSG2001}. However, to the best of the author's knowledge, the soliton solutions of the dark soliton version of the TH equations via the RH method have never been investigated by any authors.

The letter is organized as follows. In section 2, we establish a specific RH problem based on the inverse scattering transformation. In section 3, we compute N-soliton solutions of the TH equations from a specific RH problem, which possesses the identity jump matrix on the real axis. Finally, quiet a few discussions and conclusions are given in section 4.

\section{ The Riemann-Hilbert problem}
In what follows, we choose $k_1=i,A_1=\frac{i}{2}$ for the convenient of the analysis.
Then, system \eqref{1.1} possess the following Lax pair\cite{BSG2001}
\beq
\psi_x=U\psi=(\frac{i}{2}\lambda\sigma_3+iQ)\psi,\quad
\psi_t=V\psi=[-\frac{i}{2}(\epsilon\lambda^3+\lambda^2)\sigma_3+G]\psi,\label{2.1}\eeq
where $G=-i\epsilon\lambda^2Q+\lambda(i\epsilon Q^2\sigma_3-\epsilon\sigma_3Q_x-iQ)
-\sigma_3Q_x+i\epsilon Q_{xx}+i\epsilon Q^2\sigma_3+2i\epsilon Q^3+\epsilon(Q_xQ-QQ_x)$, and
\eqa\begin{array}{l}
\sigma_3=\left(\begin{array}{ccc}
-1 & 0 & 0 \\
0 & 1 & 0\\
0 & 0 & 1 \end{array} \right),\quad
Q=\left(\begin{array}{ccc}
0 & -q_1 & -q_2 \\
q_1^* & 0 & 0\\
q_2^* & 0 & 0 \end{array} \right).
\end{array}\label{2.2}\eeqa
Direct computations display that the zero-curvature equation $U_x-V_t+[U,V] = 0$ exactly gives system \eqref{1.1}.

In fact, \eqref{2.1} is equivalent to
\beq
\Phi_x-\frac{i}{2}\lambda\sigma_3\Phi=iQ\Phi,\quad
\Phi_t+\frac{i}{2}(\lambda^2+\epsilon\lambda^3)\sigma_3\Phi=G\Phi,
\label{3.1}\eeq

It is easy to see that $\tilde{A}(x,t,\lambda)=e^{\frac{i}{2}[\lambda x-(\lambda^2+\epsilon\lambda^3)t]\sigma_3}$ is a solution of the \eqref{3.1}. Introducing a new function $\Psi(x,t,\lambda)=J(x,t,\lambda)\tilde{A}(x,t,\lambda)$, and by simple calculation, we know that the spectral problems about $J(x,t,\lambda)$ satisfies
\beq
J_x-\frac{i}{2}\lambda[\sigma_3,J]=iQJ,\quad
J_t+\frac{i}{2}(\lambda^2+\epsilon\lambda^3)[\sigma_3,J]=GJ.
\label{3.3}\eeq

Now, we are construct two Jost solutions $J_{\pm}= J_{\pm}(x,\lambda)$ of the first formula in \eqref{3.3} for $\lambda\in \mathbb{R}$
\beq
J_{+}=([J_{+}]_1,[J_{+}]_2,[J_{+}]_3),\quad
J_{-}=([J_{-}]_1,[J_{-}]_2,[J_{-}]_3),
\label{3.4}\eeq
with the boundary conditions
\beq
J_{\pm}\rightarrow \mathrm{I},\,\,x\rightarrow \mp\infty.,
\label{3.5}\eeq
where $\{[J_{\pm}]_n\}_1^3$ represent the $n$-th column vector of $J_{\pm}$, $\mathrm{I}=diag\{1,1,1\}$ is a $3\times3$ unit matrix,
and the subscripts of $J(x,\lambda)$ indicate which end of the $x$-axis the boundary
conditions are set. In fact, one can using the following Volterra integral equations to define these two Jost solutions $J_{\pm}= J_{\pm}(x,\lambda)$ of the first formula in \eqref{3.3} for $\lambda\in \mathbb{R}$
\beq
J_{\pm}(x,\lambda)=\mathrm{I}+\int_{\pm\infty}^xe^{\frac{i}{2}\lambda\hat\sigma_3 (x-\xi)}iQ(\xi)J_{\pm}(\xi,\lambda)d\xi,
\label{3.6}\eeq
where $\hat\sigma_3$ is a matrix operator which acting on $3\times3$ matrix $X$ as $\hat\sigma_3 X=[\sigma_3,X]$ and
$e^{x\hat\sigma_3}X=e^{x\sigma_3}Xe^{-x\sigma_3}$.

Moreover, after simple analysis, we find that $[J_{+}]_1, [J_{-}]_2$ and $[J_{-}]_3$ admits analytic extensions to $C_{-}$. Similarly, $[J_{-}]_1, [J_{+}]_2$ and $[J_{+}]_3$ admits analytic extensions to the $C_{+}$,
here $C_{-}$ and $C_{+}$ represents the upper half $\lambda$-plane and the lower half $\lambda$-plane, respectively.

Next, one can discuss the properties of ${J}_{\pm}$. It follows from the Abel's identity and $\mathrm{Tr}(Q)=0$ that the determinants of ${J}_{\pm}$ are constants for all $x$, then from boundary conditions (2.8) yields
\beq \det {{J}_{\pm }}=1,\quad \lambda \in \mathbb{R}. \label{3.7}\eeq
In addition, we introducing another new function $ A(x,\lambda)=e^{\frac{i}{2}\lambda\sigma_3x}$, we find that spectral problem of the first formula in \eqref{3.3} exists two fundamental matrix solutions ${{J}_{+}}A$ and ${{J}_{-}}A$, which are not independent of each other but rather to enjoy a linearly correlation by a $3\times3$ scattering matrix $S(\lambda)$, that is
\beq {{J}_{-}}A={{J}_{+}}A\cdot S(\lambda),\quad \lambda \in \mathbb{R}.\label{3.8}\eeq
It follows from Eq.\eqref{3.7} and \eqref{3.8} that
\beq \det S(\lambda)=1. \label{3.10}\eeq
Moreover, let $x$ go to $+\infty$, the $3\times3$ scattering matrix $S(\lambda)$ is given as
\beq  S(\lambda)={{({{s}_{ij}})}_{3\times 3}}=\lim_{x\rightarrow+\infty}A^{-1}J_{-}A=\mathrm{I}+\int_{-\infty}^{+\infty}e^{-\frac{i}{2}\lambda\hat\sigma_3 \xi}iQJ_{-}d\xi,\,\,\lambda \in \mathbb{R}. \label{3.11}\eeq
Indeed, it follows from analytic property of $J_{-}$ that the scattering data $s_{22},s_{23},s_{32},s_{33}$ allow analytic extensions to $C_{+}$, and $s_{11}$ admits analytic extensions to $C_{-}$. Generally speaking, the other scattering data $s_{12},s_{13},s_{21}$ and $s_{31}$ cannot be extended off the
real $x$-axis.

So as to discuss behavior of Jost solutions for very large $\lambda$, we suppose
\beq J=J_0+\frac{J_1}{\lambda}+\frac{J_2}{\lambda^2}+\frac{J_3}{\lambda^3}+\frac{J_4}{\lambda^4}+\cdots \quad \lambda\rightarrow\infty,\label{3.12}\eeq
and substituting the above expansion into the first formula of \eqref{3.3} and comparing the coefficients of the same power of $\lambda$ yields
\beq\begin{array}{l}
O(\lambda^1):\frac{i}{2}[\sigma_3,J_0]=0,  \\
O(\lambda^0):J_{0,x}-\frac{i}{2}[\sigma_3,J_1]-iQJ_0=0, \\
O(\lambda^{-1}):J_{1,x}-\frac{i}{2}[\sigma_3,J_2]-iQJ_1=0,
\end{array}\label{3.13}\eeq
From $O(\lambda^1)$ and $O(\lambda^0)$ we have
\beq -\frac{i}{2}[\sigma_3,J_1]=iQJ_0,\,\,J_{0,x}=0. \label{3.14}\eeq

In order to construct the RH problem of the TH equations, we must to define another new Jost solution for the first formula of \eqref{3.3} by
\beq P_{+}=([J_{-}]_1,[J_{+}]_2,[J_{+}]_3)=J_{+}AS_{+}A^{-1}=J_{+}A\left(\begin{array}{ccc}
s_{11} & 0 & 0\\
s_{21} & 1 & 0\\
s_{31} & 0 & 1
\end{array}\right)A^{-1}, \label{3.15}\eeq
which is analytic for $\lambda\in C_{+}$ and admits asymptotic behavior for very large $\lambda$ as
\beq P_{+}\rightarrow\mathrm{I},\, \lambda\rightarrow +\infty,\, \lambda\in C_{+}. \label{3.16}\eeq

Furthermore, we also need consider the adjoint scattering
equation of the first formula \eqref{3.3}, that is
\beq \Phi_x-\frac{i}{2}\lambda[\sigma_3,\Phi]=-iQ\Phi. \label{3.17}\eeq
for the convenient of the analysis, we denote the analytic counterpart of $P_{+}$ in $C_{-}$ by $P_{-}$.
Obviously, the inverse matrices $J_{\pm}^{-1}$ defined as
\beq {[J_{+}]}^{-1}=({[J_{+}^{-1}]}^1,{[J_{+}^{-1}]}^2,{[J_{+}^{-1}]}^3)^T,\quad
{[J_{-}]}^{-1}=({[J_{-}^{-1}]}^1,{[J_{-}^{-1}]}^2,{[J_{-}^{-1}]}^3)^T,\label{3.18}\eeq
satisfy this adjoint equation \eqref{3.17}, here $[J_{\pm}^{-1}]^n(n=1,2,3)$ denote the $n$-th row vector of $J_{\pm}^{-1}$. Then we can see that ${[J_{+}^{-1}]}^1, {[J_{-}^{-1}]}^2$ and ${[J_{-}^{-1}]}^3$ admits analytic extensions to $C_{-}$. On the other hand, ${[J_{-}^{-1}]}^1, {[J_{+}^{-1}]}^2$ and ${[J_{+}^{-1}]}^3$ admits analytic extensions to the $C_{+}$.

In addition, it is not difficult to find that the inverse matrices $J_{+}^{-1}$ and $J_{-}^{-1}$ satisfy the following boundary conditions .
\beq J_{\pm}^{-1}\rightarrow \mathrm{I},\,\,x\rightarrow \mp\infty.\label{3.19}\eeq
Therefore, one can define a matrix function $P_{-}$ is expressed as follows:
\beq P_{-}=({[J_{-}^{-1}]}^1,{[J_{+}^{-1}]}^2,{[J_{+}^{-1}]}^3)^T.\label{3.20}\eeq
Through an analysis similar to the above, one can manifest that the $P_{-}$ analytic in $C_{-}$ and
\beq P_{-}\rightarrow\mathrm{I},\, \lambda\rightarrow -\infty,\, \lambda\in C_{-}. \label{3.21}\eeq
Assume that $R(k)=S^{-1}(k)$, we have
\beq {J_{-}^{-1}}=AR(\lambda)A^{-1}{J_{+}^{-1}},\label{3.22}\eeq
and
\beq P_{-}=\left(\begin{array}{ccc}
{[J_{-}^{-1}]}^1\\
{[J_{+}^{-1}]}^2\\
{[J_{+}^{-1}]}^3
\end{array}\right)=AR_{+}A^{-1}J_{+}^{-1}=A\left(\begin{array}{ccc}
r_{11} & r_{12} & r_{13}\\
0 & 1 & 0\\
0 & 0 & 1
\end{array}\right)A^{-1}J_{+}^{-1}, \label{3.23}\eeq

So far, we have established two matrix-value functions $P_{\pm}(x,\lambda)$ and $P_{-}(x,\lambda)$ which are analytic for $\lambda$ in $C_{\pm}$, respectively. In fact, these two matrix-value functions $P_{\pm}(x,k)$ can be construct a RH problem:
\beq P_{-}(x,\lambda)P_{+}(x,\lambda)=T(x,\lambda),\,\, \lambda\in C_{-}. \label{3.24}\eeq
where
\beq T(x,\lambda)=AR_{+}S_{+}A^{-1}=\left(\begin{array}{ccc}
1 & r_{12}e^{-i\lambda x} & r_{13}e^{-i\lambda x}\\
s_{21}e^{i\lambda x} & 1 & 0\\
s_{31}e^{i\lambda x} & 0 & 1
\end{array}\right),\,\, \lambda\in C_{-}. \label{3.25}\eeq
Here we have adopted the identity ${{r}_{11}}{{s}_{11}}+{{r}_{12}}{{s}_{21}}+{{r}_{13}}{{s}_{31}}=1$, and the jump contour is real $x$-axis.

Furthermore, since $J_{-}$ satisfies the temporal part of spectral equation
\eqa J_{-,t}+\frac{i}{2}(\lambda^2+\epsilon\lambda^3)[\sigma_3,J_{-}]=GJ_{-}, \label{3.26}\eeqa
we have
\eqa (\tilde{A}^{-1}J_{-}\tilde{A})_t=\tilde{A}^{-1}GJ_{-}\tilde{A},\,\,\tilde{A}=e^{i\lambda\sigma x-8i\epsilon\lambda^3\sigma_1 t} \label{3.27}\eeqa
suppose $u$ and $v$ sufficient smoothness and decay as $x\rightarrow\infty$, we have $Q_0\rightarrow 0$ as $x\rightarrow\pm\infty$.
Then taking the limit $x\rightarrow+\infty$ of Eq.\eqref{3.27} yields
\eqa S_{t}=-\frac{i}{2}(\lambda^2+\epsilon\lambda^3)[\sigma_3,S],\label{3.28}\eeqa
This above equation imply that the scattering data $s_{11},s_{22},s_{33},s_{23},s_{32}$ are time independent,
and the other scattering data satisfies
\eqa s_{1j}(t,\lambda)=s_{1j}(0,\lambda)e^{i(\lambda^2+\epsilon \lambda^3)t},\,\,
s_{j1}(t,\lambda)=s_{j1}(0,\lambda)e^{-(i\lambda^2+\epsilon \lambda^3)t},\,\,j=2,3.\label{3.30}\eeqa

\section{The soliton solutions}

From the $P_{+}$ and $P_{-}$ defined in section 2 as well as the $J_{+}$ and $J_{-}$ satisfies the scattering relationship \eqref{3.8}, it is easy to find that
\beq \mathrm{det} P_{+}(x,\lambda)=s_{11}(\lambda),\quad \mathrm{det} P_{-}(x,\lambda)=r_{11}(\lambda),\label{4.1}\eeq
where $r_{11}=s_{22}s_{33}-s_{23}s_{32}$. Indeed, owing to the $s_{11}$ and $r_{11}$ are time independent, then the zeros of $s_{11}=0$ and $r_{11}=0$ are also time independent. Moreover, due to $\sigma_3Q\sigma_3=-Q$ and $Q^\dag=-Q$, it is not difficulty to find that
\beq J^\dag(x,t,\lambda^*)=\sigma_3 J^{-1}(x,t,\lambda)\sigma_3,\quad S^\dag(\lambda)=\sigma_3 S^{-1}(\lambda)\sigma_3,\label{4.2}\eeq
then
\beq P_{+}^{\dag}(\lambda)=\sigma_3P_{-}\sigma_3,\label{4.3}\eeq

Suppose that $s_{11}$ possess $N\geq0$ possible zeros in $C_{+}$ denoted by $\{\lambda_j,1\leq j\leq N\}$, and $r_{11}$ possess $N\geq0$ possible zeros in $C_{-}$ denoted by $\{\hat\lambda_j,1\leq j\leq N\}$. For the sake of simplicity, one can suppose that all zeros$\{(\lambda_j,\hat \lambda_j),j=1,2,...,N\}$ of $s_{11}$ and $r_{11}$ are simple zeros. In this case, each of kernel $P_{+}(\lambda_j)$ and kernel $P_{-}(\hat \lambda_j)$ include only a single column vector $v_{j}$ and row vector $\hat v_{j}$, respectively, such that
\eqa P_{+}(\lambda_j)v_{j}=0,\quad \hat v_{j}P_{-}(\hat \lambda_j)=0.  \label{4.4}\eeqa

Owing to $P_{+}(\lambda)$ is the solution of the first formula of \eqref{3.3}, we assume that the asymptotic expansion of $P_{+}(\lambda)$ at large $\lambda$ as
\beq P_{+}=\mathrm{I}+\frac{P_{+}^{(1)}}{\lambda}+O(\lambda^{-2})\quad \lambda\rightarrow\infty,\label{4.5}\eeq
substitute the above expansion into \eqref{3.3} and compare $O(1)$ terms obtain
the potential functions $q_1$ and $q_2$ can be reconstructed by
\beq q_1=-(P_{+}^{(1)})_{12},\quad q_2=-(P_{+}^{(1)})_{13},\label{4.6}\eeq
where $P_{+}^{(1)}=(P_{+}^{(1)})_{3\times3}$ and $(P_{+}^{(1)})_{ij}$ is the $(i;j)$-entry of $P_{+}^{(1)},i,j=1,2,3$.

In order to obtain the spatial evolutions for vectors $v_j(x,t)$, on the one hand, we taking the $x$-derivative to equation $P_{+}v_j=0$ and using the first formula of \eqref{3.3} obtain
\beq v_{kx}=\frac{i}{2}\lambda_k\sigma_3v_k,\label{4.7}\eeq
on the other hand, we also taking the $t$-derivative to equation $P_{+}v_j=0$ and using the second formula of \eqref{3.3} obtain
\beq v_{kt}=-\frac{i}{2}(\lambda_k^2+\epsilon\lambda_k^3)\sigma_3v_k,\label{4.8}\eeq
By solving \eqref{4.7} and \eqref{4.8} explicitly, we get
\beq v_{k}(x,t)=e^{\frac{i}{2}[\lambda_kx-(\lambda_k^2+\epsilon\lambda_k^3)t]\sigma_3}v_{k0},\,\,
 \hat v_{k}(x,t)=v_{k0}^\dag e^{\frac{i}{2}(-\lambda_k^*x+[(\lambda_k^*)^2+\epsilon(\lambda_k^*)^3]t)\sigma_3}\sigma_3.\label{4.9}\eeq
where $v_{k0}$ and $\hat v_{k0}$ are constant vectors.

In order to obtain multi-soliton solutions for the TH equations \eqref{1.1}, one can choose the jump matrix $T=\mathrm{I}$ is a $3\times3$ unit matrix in \eqref{3.24}. That is to say, the discrete scattering data $r_{12}=r_{13}=s_{21}=s_{31}=0$, consequently, the unique solution to this special RH problem have been solved in \cite{YJK2010}, and the result is
\beq
P_{+}(k)=\mathrm{I}-\sum_{j=1}^{N}\sum_{k=1}^{N}\frac{v_j(M^{-1})_{jk}\hat v_k}{\lambda-\hat \lambda_j}.
\label{4.10}\eeq
where $M=(M_{jk})_{N\times N}$ is a matrix whose entries are
\eqa M_{jk}=\frac{\hat v_jv_k}{\lambda_j^*-\lambda_k},\,1\leq j,k\leq N.\label{4.11}\eeqa
Therefore, from \eqref{4.10}, we obtain
\beq
P_{+}^{(1)}=\sum_{j=1}^{N}\sum_{k=1}^{N}v_j(M^{-1})_{jk}\hat v_k.\label{4.12}
\eeq
It follows from \eqref{4.12} that the general N-soliton solution for the TH equations \eqref{1.1} reads
\beq
q_1=-\sum_{j=1}^{N}\sum_{k=1}^{N}m_k^*e^{-\theta_j+\theta_k^*}(M^{-1})_{jk},\,\,
q_2=-\sum_{j=1}^{N}\sum_{k=1}^{N}n_k^*e^{-\theta_j+\theta_k^*}(M^{-1})_{jk}.\label{4.13}
\eeq
and $M=(M_{jk})_{N\times N}$ is given by
\eqa M_{jk}=\frac{-e^{-(\theta_j^*+\theta_k)}+(m_j^*m_k+n_j^*n_k)e^{\theta_j^*+\theta_k}}{\lambda_j^*-\lambda_k},\,1\leq j,k\leq N.\label{4.14}\eeqa
with $\theta_k=\frac{i}{2}[\lambda_kx-(\lambda_k^2+\epsilon\lambda_k^3)t]$,
we have chosen $v_{k0}=[1,m_k,n_k]^T$.

Then, on the one hand, as a special example, one can choose $N=1$ in formula \eqref{4.13} and with \eqref{4.11}, we obtain the one-soliton solution as follows:
\beq
q_1(x,t)=-\frac{m_1^*e^{\theta_1^*-\theta_1}(\lambda_1^*-\lambda_1)}{-e^{-(\theta_1^*+\theta_1)}+(|m_1|^2+|n_1|^2)e^{\theta_1^*+\theta_1}},\,\,
q_2(x,t)=-\frac{n_1^*e^{\theta_1^*-\theta_1}(\lambda_1^*-\lambda_1)}{-e^{-(\theta_1^*+\theta_1)}+(|m_1|^2+|n_1|^2)e^{\theta_1^*+\theta_1}}.
\label{4.15}\eeq
Letting $\lambda_1=\lambda_{11}+i\lambda_{12}$, then the one-soliton solution \eqref{4.15} can be written as
\beq
q_1(x,t)=-i\lambda_{12}m_1^*e^{\theta_1^*-\theta_1-\xi_1}\mathrm{csch}(\theta_1^*+\theta_1+\xi_1),\,\,\,
q_2(x,t)=-i\lambda_{12}n_1^*e^{\theta_1^*-\theta_1-\xi_1}\mathrm{csch}(\theta_1^*+\theta_1+\xi_1).
\label{4.16}\eeq
where $\theta_1^*-\theta_1=-i[\lambda_{11}x-(\lambda_{11}^2-\lambda_{12}^2)t-\epsilon(\lambda_{11}^3-3\lambda_{11}\lambda_{12}^2)t]$,
$\theta_1^*+\theta_1=-\lambda_{12}x+2\lambda_{11}\lambda_{12}t+\epsilon(3\lambda_{11}^2\lambda_{12}-\lambda_{12}^3)t$ and $\xi_1$ satisfy $e^{2\xi_1}=|m|_1^2+|n|_1^2$.

On the other hand, as another special example, one can choose $N=2$ in formula \eqref{4.13} and with \eqref{4.11}, we arrive at the two-soliton solution as follows:
\beq \begin{array}{l}
q_1(x,t)=-[m_1^*e^{\theta_1^*-\theta_1}(M^{-1})_{11}+m_1^*e^{\theta_1^*-\theta_2}(M^{-1})_{21}
+m_2^*e^{\theta_2^*-\theta_1}(M^{-1})_{12}+m_2^*e^{\theta_2^*-\theta_2}(M^{-1})_{22}],\\
q_2(x,t)=-[n_1^*e^{\theta_1^*-\theta_1}(M^{-1})_{11}+n_1^*e^{\theta_1^*-\theta_2}(M^{-1})_{21}
+n_2^*e^{\theta_2^*-\theta_1}(M^{-1})_{12}+n_2^*e^{\theta_2^*-\theta_2}(M^{-1})_{22}].
\end{array}\label{4.17}\eeq
where $\theta_k=\frac{i}{2}[\lambda_kx-(\lambda_k^2+\epsilon\lambda_k^3)t]$, $M=(M_{jk})_{2\times 2}$ with
\beq \begin{array}{l}
M_{11}=\frac{2e^{\xi_1}}{\lambda_1^*-\lambda_1}\mathrm{sinh}(\theta_1^*+\theta_1+\xi_1),\,\,
M_{12}=\frac{2e^{\xi_2}}{\lambda_1^*-\lambda_2}\mathrm{sinh}(\theta_1^*+\theta_2+\xi_2),\\
M_{21}=\frac{2e^{\xi_2^*}}{\lambda_2^*-\lambda_1}\mathrm{sinh}(\theta_1+\theta_2^*+\xi_2^*),\,\,
M_{22}=\frac{2e^{\xi_3}}{\lambda_2^*-\lambda_2}\mathrm{sinh}(\theta_2^*+\theta_2+\xi_3),
\end{array}\nn\eeq
and $e^{2\xi_1}=m_1^*m_2+n_1^*n_2$ and $e^{2\xi_3}=|m|_2^2+|n|_2^2$, respectively.

\section{Discussions and conclusions}

In fact, as a promotion, the integrable two-component Hirota equations can be extended to the integrable generalized multi-component Hirota equations:
\beq k_1q_{lt}+2A_1k_1q_{lxx}+4k_1^3A_1\sum_{j=1}^N|q_j|^2q_l\\
-i\epsilon[-k_1q_{lxxx}-3ik_1^3\sum_{j=1}^N|q_j|^2q_{lx}-3ik_1^3q_l\sum_{j=1}^Nq_j^*q_{jx}]=0,l=1,2,\ldots,N
\eeq
which possess the following Lax pair
\beq
\psi_x=U\psi=(\frac{i}{2}\lambda\sigma_3+iQ)\psi,\quad
\psi_t=V\psi=[\frac{i}{2}(\epsilon\lambda^3+\lambda^2)\sigma_3+G]\psi,\label{5.2}\eeq
where $G=-i\epsilon\lambda^2Q+\lambda(i\epsilon Q^2\sigma_3-\epsilon\sigma_3Q_x-iQ)
-\sigma_3Q_x+i\epsilon Q_{xx}+i\epsilon Q^2\sigma_3+2i\epsilon Q^3+\epsilon(Q_xQ-QQ_x)$, and
\eqa\begin{array}{l}
\sigma_3=\left(\begin{array}{cc}
-1&\mathbf{0}_{1\times N}\\
\mathbf{0}_{N\times 1}&\mathbf{I}_{N\times N}\end{array} \right),
Q=\left(\begin{array}{ccc}
0&-\mathbf{q}^T\\
\mathbf{q}^*&\mathbf{0}_{N\times N}\end{array} \right).
\end{array}\label{5.3}\eeqa
with $\mathbf{q}=(q_1,q_2,...,q_N)^T$.
When $k_1=1,A_1=-\frac{i}{2}$, which means to the bright soliton version of the multi-component Hirota equations, when $k_1=i,A_1=\frac{i}{2}$, which means to the dark soliton version of the multi-component Hirota equations. Accordingly, one can also examine the N-soliton solutions to the integrable generalized multi-component Hirota equations by the same way in above two Section. However, we don't examine them here since the procedure is mechanical.

Moreover, based on the $3\times 3$ matrix RH problem of the TH equations discussed by authors in \cite{Hu2}, one can examine the long-time asymptotic behavior for the solutions of the TH equations via the nonlinear steepest descent method introduce by Deift and Zhou \cite{Deift1993}.

\section*{Acknowledgements}

The work was supported by the NSF of China under Grant Nos.11601055, 11805114, NSF of Anhui Province under Grant No.1408085QA06.

\end{document}